\documentclass{elsarticle}
\usepackage{srcltx,
enumerate,
amsmath,amssymb}
\input{xy}
\xyoption{all}

\usepackage{amsthm,
varioref,
graphicx}
\usepackage[latin1]{inputenc}

\usepackage{hyperref}

		\newcommand{\category}[1]{\mathbf{#1}}
	
    \newcommand{\D}{\category {D}} 

   \newcommand{\gm}{\mathrm{gm}}

		\newcommand{\Shvv}{\category{Shv}} 

    \newcommand{\DM}{\category {DM}}

\newcommand{\DMgm}{\DM_\gm} 

    \newcommand{\DTM}{\category{DTM}} 
    \newcommand{\MTM}{\category{MTM}} 
    \newcommand{\DATM}{\category{DATM}} 
    \newcommand{\MATM}{\category{MATM}} 

    \newcommand{\MM}{\category {MM}} 

	    \newtheoremstyle{Normal}
  {}
  {}
  {}
  {}
  {\bfseries}
  {.}
  { }
  {}

    \theoremstyle{Normal} 

    \newtheorem{Defi}{Definition}[section]

    \newtheorem{Conj}[Defi]{Conjecture}
    \newtheorem{Bem}[Defi]{Remark}
    \newtheorem{Bsp}[Defi]{Example}
     
    \newtheorem{Axio}[Defi]{Axiom}
    \newtheorem{Ques}[Defi]{Question}

		\theoremstyle{remark}
	
    \theoremstyle{plain} 

    \newtheorem{Satz}[Defi]{Proposition}
    \newtheorem{Theo}[Defi]{Theorem}
    \newtheorem{Folg}[Defi]{Corollary}
    \newtheorem{Lemm}[Defi]{Lemma}

\newcommand{\refit}[1]{(\ref{item_#1})}

\newcommand{\refsect}[1]{Section \ref{sect_#1}}
\newcommand{\conj}{\begin{Conj}} 			\newcommand{\xconj}{\end{Conj}}											    
\newcommand{\ques}{\begin{Ques}} 			\newcommand{\xques}{\end{Ques}}											    
\newcommand{\axio}{\begin{Axio}} 			\newcommand{\xaxio}{\end{Axio}}											    
\newcommand{\bem}{\begin{Bem}} 			\newcommand{\xbem}{\end{Bem}}											    \newcommand{\refbe}[1]{Remark \ref{bem_#1}}
\newcommand{\defi}{\begin{Defi}} 			\newcommand{\xdefi}{\end{Defi}}										\newcommand{\refde}[1]{Definition \ref{defi_#1}}
\newcommand{\lemm}{\begin{Lemm}}			\newcommand{\xlemm}{\end{Lemm}}											\newcommand{\refle}[1]{Lemma \ref{lemm_#1}}
\newcommand{\satz}{\begin{Satz}}			\newcommand{\xsatz}{\end{Satz}}										\newcommand{\refsa}[1]{Proposition \ref{satz_#1}}
\newcommand{\theo}{\begin{Theo}}			\newcommand{\xtheo}{\end{Theo}}											\newcommand{\refth}[1]{Theorem \ref{theo_#1}}
\newcommand{\bsp}{\begin{Bsp}}				\newcommand{\xbsp}{\end{Bsp}}												
\newcommand{\folg}{\begin{Folg}}				\newcommand{\xfolg}{\end{Folg}}
\newcommand{\coro}{\folg}				\newcommand{\xcoro}{\xfolg} 			\newcommand{\refcor}[1]{Corollary \ref{coro_#1}}
\newcommand{\cor}{\begin{Folg}}				\newcommand{\xcor}{\end{Folg}}
\newcommand{\mycomment}{}

\newcommand{\eqnarra}{\begin{eqnarray}}				\newcommand{\xeqnarra}{\end{eqnarray}}
\newcommand{\eqnarr}{\begin{eqnarray*}}				\newcommand{\xeqnarr}{\end{eqnarray*}}

\newcommand{\eqn}{\begin{equation}} 		\newcommand{\xeqn}{\end{equation}}
\newcommand{\refeq}[1]{(\ref{eqn_#1})}

\newcommand{\gl}{\index}

\newcommand{\fcohoSectionThree}{Section 3}
\newcommand{\fcohoSectionFour}{Section 4}

  \newcommand{\Dph}[1]{\emph{#1}}
  \newcommand{\Def}[1]{\Dph{#1}\gl{#1}}

\newcommand{\Demm}[2]{$#2$\index{#1@$#2$}}

\newcommand{\mylabel}[1]{\label{#1}}

\newcommand{\status}[1]{}

\newcommand{\OF}{{\mathcal{O}_F}}

\newcommand{\R}{\mathrm{R}}
\newcommand{\RG}{{\R} {\Gamma}}

\newcommand{\lr}{{\longrightarrow}}
\renewcommand{\r}{\rightarrow}

\renewcommand{\t}{{\otimes}}

\newcommand{\A}{\mathbb{A}^1}
\renewcommand{\P}[1][1]{\mathbb P^{#1}}

\newcommand{\Q}{\mathbb{Q}}

\newcommand{\Ql}{{\mathbb{Q}_{{\ell}}}}

\newcommand{\Z}{\mathbb{Z}}

\newcommand{\pp}{\mathfrak{p}}
\newcommand{\Fpp}{{\mathbb{F}_{\pp}}}

\newcommand{\Fq}{{\mathbb{F}_q}}

\renewcommand{\H}{\mathrm{H}}

\newcommand{\p}{{{}^\mathrm{p}}}
\newcommand{\pH}{{\p\H}}

\newcommand{\x}{{\times}}
\newcommand{\ol}[1]{{\overline{#1}}} 
\newcommand{\one}{\mathbf{1}}

\newcommand{\Beweis}{{\normalfont} \textbf{Proof}}

\newcommand{\lc}{\textit{loc.~cit.}}
\newcommand{\oc}{\textit{op.~cit.}}
\newcommand{\lcs}{\textit{loc.~cit.}\ }

\newcommand{\Spec}{\mathrm{Spec}\text{ }}
\newcommand{\SpecOF}{{\Spec} {\OF}}

\newcommand{\SpecZ}{{\Spec} {\Z}}

\newcommand{\SpecFq}{{\Spec} {\Fq}}

\newcommand{\SpecF}{{\Spec} {F}}

\newcommand{\SpecFpp}{{\Spec} {\Fpp}}

\newcommand{\op}{\mathrm{op}}

\newcommand{\eff}{\mathrm{eff}}

\newcommand{\gr}{\operatorname{gr}} 

    \def\acute{\mathaccent"7013 }
    \newcommand{\et}{\mathrm{\acute et}}

\newcommand{\h}{\mathrm{h}}

\newcommand{\Hom}{\mathrm{Hom}}
\newcommand{\End}{\mathrm{End}}
\newcommand{\IHom}{\underline{\Hom}}
\newcommand{\id}{\mathrm{id}}
\renewcommand{\gr}{\operatorname{gr}} 

\newcommand{\im}{\operatorname{im}} 

\newcommand{\M}{\operatorname{M}}

\newcommand{\Gal}{\mathrm{Gal}}

\renewcommand{\h}{\mathrm{h}}
\renewcommand{\mod}[1]{\ \ (\mathrm{mod \ }#1)}

\newcommand{\bound}{\mathrm{b}} 
\newcommand{\coker}{\operatorname{coker}}
\newcommand{\Ext}{\operatorname{Ext}}

\newcommand{\cone}{\operatorname{cone}}

\newcommand{\loc}[2]{[}

\newcommand{\pr}{\begin{proof}[\Beweis: ]}
\newcommand{\pf}{\pr}
\newcommand{\xpf}{\end{proof}}
\newcommand{\xproof}{\end{proof}}

\newcommand{\wei}[1]{\langle #1 \rangle}
\newcommand{\Tw}[1]{T_{\wei{#1}}}
\newcommand{\Extcl}{\overline \Ext}

\makeindex

\makeatletter
\@addtoreset{Defi}{chapter}
\@addtoreset{Defi}{section}
\makeatother

\begin{document}

\begin{frontmatter}

\title{Mixed Artin-Tate motives over number rings}
\author{Jakob Scholbach}
\address{Universit{\"a}t M{\"u}nster, Mathematisches Institut, Einsteinstr. 62, D-48149 M{\"u}nster, Germany, 
jakob.scholbach@uni-muenster.de
}

\begin{abstract}
This paper studies Artin-Tate motives over bases $S \subset \SpecOF$, for a number field $F$. As a subcategory of motives over $S$, the triangulated category of Artin-Tate motives $\DATM(S)$ is generated by motives $\phi_* \one(n)$, where $\phi$ is any finite map. After establishing stability of these subcategories under pullback and pushforward along open and closed immersions, a motivic $t$-structure is constructed. Exactness properties of these functors familiar from perverse sheaves are shown to hold in this context. The cohomological dimension of mixed Artin-Tate motives ($\MATM(S)$) is two, and there is an equivalence $\DATM(S) \cong \D^\bound(\MATM(S))$. Finally, mixed Artin-Tate motives enjoya strict functorial weight filtration.
\end{abstract}


\end{frontmatter}

Geometric motives, as developed by Hanamura, Levine, and Voevodsky \cite{Hanamura:Motives1, Levine:MixedMotives, Voevodsky:TCM}, are established as a valuable tool in understanding geometric and arithmetic aspects of algebraic varieties over fields. However, the stupefying ambiance inherent to motives, exemplified by Grothendieck's motivic proof idea of the Weil conjectures, remains largely conjectural---especially what concerns the existence of mixed motives $\MM(K)$ over some field $K$. That category should be the heart of the so-called motivic $t$-structure on $\DMgm(K)$, the category of geometric motives. Much the same way as the cohomology groups of a variety $X$ over $K$, e.g.\ $\H^n_\et(X \x_K \ol K, \Ql)$, $\ell$-adic cohomology for $\ell \neq \operatorname{char} K$ are commonly realized as cohomology groups of a complex, e.g.\ $\RG_\ell(X, \Ql)$, there should be mixed motives $\h^n(X)$ that are obtained by applying truncation functors belonging to the $t$-structure to $\M(X)$, the motive of $X$. However, progress on mixed motives has proved hard to come by. To date, such a formalism has been developed for motives of zero- and one-dimensional varieties, only. This is due to Levine, Voevodsky, Orgogozo and Wildeshaus \cite{Levine:TateMotives, Voevodsky:TCM, Orgogozo:Isomotifs, Wildeshaus:ATM}.

Building upon Voevodsky's work, Ivorra and recently Cisinski and D\'eglise \cite{Ivorra:These, CisinskiDeglise:Triangulated} developed a theory of geometric motives $\DMgm(S)$ over more general bases. The purpose of this work is to join the ideas of Beilinson, Bernstein and Deligne on perverse sheaves \cite{BBD} with the ones on Artin-Tate motives over fields to obtain a workable category of mixed Tate and Artin-Tate motives over bases $S$ which are open subschemes of $\SpecOF$, the ring of integers in a number field $F$. As over a field, this provides some piece of evidence for the existence and properties of the conjectural category of mixed motives over $S$ and its properties.

The triangulated category $\DTM(S)$ ($\DATM(S)$) of Tate (Artin-Tate) motives is defined (\ref{defi_DATM}) to be the triangulated subcategory of $\DMgm(S)$ (with rational coefficients) generated by direct summands of $\one(n)$ and $i_*\one(n)$ ($\phi_* \one(n)$, respectively). Here, $\one$ is a shorthand for the motive of the base scheme, $(n)$ denotes the Tate twist, $i: \SpecFpp \r S$ is a closed point, $\phi: V \r S$ is any finite map and $\phi_*$ etc.\ denotes the pushforward functor on geometric motives $\DMgm(V) \r \DMgm(S)$. In case $S$ is a finite disjoint union of $\SpecFpp$, the usual definition of (Artin\nobreakdash-)Tate motives over $S$ is recalled in \refde{ATfields}.   

The following theorem and its ``proof'' is an overview of the paper. 

\theo
The categories $\DTM(S)$ and $\DATM(S)$ are stable under standard functoriality operations such as $i^!$, $j_*$ etc.\ for open and closed embeddings $j$ and $i$, respectively. 

Both categories enjoy a non-degenerate $t$-structure called \emph{motivic $t$-structure}. Its heart is denoted $\MTM(S)$ or $\MATM(S)$, respectively and called category of \emph{mixed (Artin\nobreakdash-)Tate motives}. 

The functors $i^*$, $j_*$ etc.\ feature exactness properties familiar from the corresponding situation of perverse sheaves. For example, $i^!$ is left-exact, and $j_*$ is exact with respect to the motivic $t$-structure. 

The cohomological dimension of $\MTM(S)$ and $\MATM(S)$ is one and two, respectively. We have an equivalence of categories
$$\D^\bound(\MATM(S)) \cong \DATM(S)$$
and likewise for Tate motives.

The ``site'' of mixed Artin-Tate motives over $S$ has enough points in the sense that a mixed Artin-Tate motive over $S$ is zero if and only if its restrictions to all  closed points of $S$ vanish. 

The non-unique weight truncation triangles on $\DATM(S)$ \`a la Bondarko can be refined to a strict functorial weight filtration on $\MATM(S)$.  
\xtheo

\pf
The first statement is \refle{respect}. It is proven using the localization, purity and base-change properties of geometric motives.

We will write $T(S)$ for either $\DTM(S)$ or $\DATM(S)$. The existence of the motivic $t$-structure on $T(S)$ is proven in three steps. The first ingredient is the well-known motivic $t$-structure on Artin-Tate motives over finite fields (\refle{tstructureFp}). The second step is the study of a subcategory $\tilde T(S) \subset T(S)$ generated by $\phi_* \one(n)$, where $\phi$ is finite and \'etale (Artin-Tate motives), or just by $\one(n)$ (Tate motives). This category is first equipped with an an auxiliary $t$-structure. Then, a motivic $t$-structure on $\tilde T(S)$ is defined in \refsect{motivict} by using the cohomology functor for the auxiliary $t$-structure. This statement uses (and its proof imitates) the corresponding situation for Artin-Tate motives over number fields due to Levine and Wildeshaus. The $t$-structure on $\tilde T(S)$ is glued with the one over finite fields, using the general gluing procedure of $t$-structures of \cite{BBD}, see \refth{motivictstructure}. Much the same way as with perverse sheaves, there are shifts accounting for $\dim S = 1$, that is to say, $i_* \one(n)$ and $\one(n)[1]$ are mixed Tate motives. Beyond the formalism of geometric motives, the only non-formal ingredient of the motivic $t$-structure are vanishing properties of algebraic $K$-theory of number rings, number fields and finite fields due to Quillen, Borel and Soul\'e.

The exactness statements are shown in \refth{exactness}. This theorem gives some content to the exactness axioms for general mixed motives over $S$ \cite[\fcohoSectionFour]{Scholbach:fcoho}. The key stepstone is the following: for any immersion of a closed point $i: \SpecFpp \r S$, the functor $i^*$ maps the heart $T^0(S)$ of $T(S)$ to $T^{[-1,0]}(\SpecFpp)$, that is, the category of (Artin\nobreakdash-)Tate motives over $\Fpp$ whose only nonzero cohomology terms are in degrees $-1$ and $0$. The proof is a careful reduction to basic calculations which relies on facts gathered in \refsect{motivict} about the heart of $\tilde T(S)$. 

The cohomological dimensions are calculated in \refsa{cohomdimAT}. The Artin-Tate case is a special (but non-conjectural) case of a similar fact for general mixed motives over $S$. The difference in the Tate case is because the generators of $\DTM(S)$ have good reduction at all places. 

By an argument of Wildeshaus, under a mild homological condition on $T(S)$, the identity on $T^0(S)$ extends to a functor $\D^\bound(T^0(S)) \r T(S)$ (\refth{Wildeshaus}). While it is an equivalence in the case of Tate motives for formal reasons, the Artin-Tate case requires some localization arguments.

The last but one statement is \refsa{enoughpoints}. It might be seen as a first step into motivic sheaves. 

The weight filtration is established in \refsect{weights}. The key idea is an in-depth analysis of a particular family generators, namely motives of the form $f_* \one[1]$, where $f$ is a finite map between regular schemes. 
\xpf

Deligne and Goncharov define a category of mixed Tate motives over rings $\mathcal O_S$ of $S$-integers of a number field $F$ \cite[1.4., 1.7.]{DeligneGoncharov}. Unlike the mixed Tate motives we study, their category is a \emph{sub}category of mixed Tate motives over $F$, consisting of motives subject to certain non-ramification constraints, akin to Scholl's notion of mixed motives over $\OF$ \cite{Scholl:Remarks}.

This paper is an outgrowth of part of my thesis. I owe many thanks to Annette Huber for her advice during that time. I am also grateful to Denis-Charles Cisinski and Fr\'ed\'eric D\'eglise for teaching me their work on motives over general bases and to the referee for suggesting that \refsect{weights} be (re-)written.


\section{Geometric motives} \mylabel{sect_preliminaries}
In this section we briefly recall some properties of the triangulated categories of geometric motives $\DMgm(X)$, where $X$ will be either a number field $F$ or an open or closed subscheme of $\SpecOF$. All of this is due to Cisinski and D\'eglise \cite{CisinskiDeglise:Triangulated}. 

The categories $\DMgm(X)$, where $X$ is any of the afore-mentioned bases, are related by adjoint functors 
$f^*: \DMgm(X) \leftrightarrows \DMgm(Y) : f_*$, where $f: Y \r X$ is any map, and 
$f_!: \DMgm(Y) \leftrightarrows \DMgm(X) : f^!$ ($f$ separated of finite type). If $f$ is smooth, $f^*$ also has a left adjoint $f_\sharp$. The category $\DMgm(X)$ enjoys inner $\Hom$'s, denoted $\IHom$, and a tensor structure whose unit is denoted $\one$. Pullback functors $f^*$ are monoidal. In particular $f^* \one_X = \one_Y$ for $f: Y \r X$. The \Def{motive} of any scheme $f: Y \r X$ of finite type is defined as $f_! f^! \one$ and denoted $\M(Y)$. (For $f$ smooth, \cite[Section 1.1.]{CisinskiDeglise:Triangulated} puts $\M(Y) := f_\sharp f^* \one$. This agrees with the previous definition by relative purity, see below.) The tensor structure in $\DMgm(X)$ is such that 
\eqn
\mylabel{eqn_monoidal}
\M(Y) \t \M(Y') = \M(Y \x_X Y')
\xeqn 
for any two smooth schemes $Y$ and $Y'$ over $X$. There is a distinguished object $\one(1)$ such that $\M(\P_X) = \one \oplus \one(1)[2]$. Tensoring with $\one(1)$ is an equivalence on $\DMgm(X)$, and $\one(n)$ is defined in the usual way in terms of tensor powers of $\one(1)$. We exclusively work with rational coefficients, i.e., all morphism groups are $\Q$-vector spaces. If $X$ is regular, morphisms in $\DMgm(X)$ are given by
\eqn
\mylabel{eqn_morphismsKTheory}
\Hom_{\DMgm(X)}(\one, \one(q)[p]) \cong K_{2q-p}(X)_\Q^{(q)},
\xeqn 
the $q$-th Adams eigenspace in algebraic $K$-theory of $X$, tensored with $\Q$ \cite[Section 13.2]{CisinskiDeglise:Triangulated}. Having rational coefficients (or coefficients in a bigger number field) is vital when it comes to vanishing properties of $\Hom$-groups in $\DMgm(X)$. (With integral coefficients, the existence of a $t$-structure even in the case of Artin motives over a field is unclear.) 

Throughout we need a property called \Def{localization}: for any closed immersion $i: Z \r X$ with open complement $j$ we have the following functorial distinguished triangles in $\DMgm(X)$
\eqn
\mylabel{eqn_localization}
j_! j^* \r \id \r i_* i^*
\xeqn

We need to know that the functors $f_!$ and $f_*$ naturally agree for any proper map $f$, as do $f^!$ and $f^*(d)[2d]$ when $f$ is smooth and quasi-projective of constant relative dimension $d$ (\Def{relative purity}). Moreover, when $i: Z \r X$ is a closed immersion of constant relative codimension $c$ and $Z$ and $X$ are regular, we have $i^! \one \cong i^* \one (-c)[-2c]$. This is called \Def{absolute purity} \cite[Sections 2.4, 13.4]{CisinskiDeglise:Triangulated}. Finally, for $f: Y \r X$, $g: X' \r X$, $f' : Y' := X' \x_X Y \r X'$ and $g': Y' \r Y$, there is a natural \Def{base-change} isomorphism of functors $f^* g_! \cong g'_! f'^*$ \cite[Section 2.2]{CisinskiDeglise:Triangulated}, originally due to Ayoub \cite{Ayoub:Six1}.

The \Def{Verdier dual} functor $D_X: \DMgm(X)^\op \r \DMgm(X)$ is defined by $D_X(M) := \IHom (M, \pi^! \one(1)[2])$ for any $M \in \DMgm(X)$, where $\pi: X \r \SpecZ$ denotes the structural map. For example, for an open subscheme $X$ of $\SpecOF$ the factorization 
$$X \subset \SpecOF \r \mathbb A^n_\Z \r \SpecZ$$ 
and absolute and relative purity show that $D_X(-) = \IHom (-, \one(1)[2])$. For $X = \SpecFq$ one gets $D_X(-) = \IHom (M, \one)$. The Verdier dual functor exchanges ``$!$'' and ``$*$'', e.g., there are natural isomorphisms $D (f^! M) \cong f^* D(M)$ \cite[Section 14.3]{CisinskiDeglise:Triangulated}. For example, the Verdier dual of \refeq{localization} yields a distinguished triangle
\eqn
\mylabel{eqn_localization2}
i_* i^! \r \id \r j_* j^*.
\xeqn
For $X = \SpecOF$, taking the limit over increasingly small open subschemes, one obtains a distinguished triangle in $\DM(X)$ of the following form \cite[Section 14.2]{CisinskiDeglise:Triangulated}. (The category $\DM(X)$ is a bigger category whose subcategory of compact objects is $\DMgm(X)$.)
\eqn
\mylabel{eqn_localizationgeneric}
\oplus_{\pp \in S} {i_\pp}_* i_\pp^! \r \id \r \eta_* \eta^*,
\xeqn
where $\eta: \SpecF \r \SpecOF$ is the generic point, the sum runs over all closed points $\pp \in X$, $i_\pp$ is the closed immersion.  

\section{Triangulated Artin-Tate motives}\mylabel{sect_definitionstability}


Recall the following classical definition. We apply it to a number field or a finite field: 
\defi \mylabel{defi_ATfields}
Let $K$ be a field. The category of Tate motives $\DTM(K)$ over $K$ is by definition the triangulated subcategory of $\DMgm(K)$ generated by $\one(n)$ where $n \in \Z$. The smallest full triangulated subcategory $\DATM(K)$ stable under tensoring with $\one(n)$ and containing direct summands of motives $f_* \one$, where $f: K' \r K$ is any finite map, is called category of Artin-Tate motives over $K$. For a scheme $S$ of the form $S = \sqcup \Spec K_i$, a finite disjoint union of spectra of fields, we put $\DATM(S) := \oplus_i \DATM(K_i)$ and likewise for $\DTM$.  
\xdefi

This section gives a generalization of that definition to bases $S$ which are open subschemes of $\SpecOF$ based on the idea that Artin-Tate motives over $S$ should be compatible with the ones over $F$ and $\Fpp$ under standard functoriality. 

\defi \mylabel{defi_DATM}
The categories $\DTM(S) \subset \DMgm(S)$ of \Def{Tate motives} and $\DATM(S) \subset \DMgm(S)$\index{DATM@$\DATM(-)$} of \Def{Artin-Tate motives} over $S$ 
are the triangulated subcategories generated by the direct summands of  
$$\one(n),\ i_*\one(n) \ \ \text{(Tate motives)}$$
and 
$$\phi_* \one(n), \ \ \text{(Artin-Tate motives)}$$ 
respectively, where $n \in \Z$, $\phi : V \r S$ is any finite map (including those that factor over a closed point) and $i$ is the immersion of any closed point of $S$.
\xdefi

\bem \mylabel{bem_onepoint}
\begin{itemize}
\item 
In comparison to motives over a field, a category of Artin motives over $S$, defined by removing the twists, is less viable, since it is not stable under Verdier duality and $i_* i^!$, where $i$ is a closed embedding of a point.
\item
We can assume by localization (see \refeq{localization}, \refeq{localization2}) that the domain of $\phi$ is a reduced scheme. 
\item The category of Tate motives $\DTM(S)$ agrees with the triangulated category generated by the above generators (without taking direct summands), see \refle{summands}.
\end{itemize}
\xbem

\noindent \emph{For brevity, we write $T(S)$ or $T$ for $\DATM(S)$ or $\DTM(S)$ in the sequel. In most proofs, we will only spell out the case of Artin-Tate motives.}\\


\lemm \mylabel{lemm_respect} 
Let $j: S' \r S$ be any open immersion, $i: Z \r S$ be any closed immersion and $f: V \r S$ any finite map such that $V$ is regular. Let $\eta: \SpecF \r S$ be the generic point. Then the functors $f_* = f_!$, $f^*$ and $f^!$ preserve Artin-Tate motives. Similar statements hold for Artin-Tate and Tate motives for $j$ and $i$. Moreover, $\eta^*$, the Verdier dual functor $D$ and the tensor product on $\DMgm(S)$ respect the subcategories of (Artin\nobreakdash-)Tate motives.
\xlemm

The functor $\eta_*$ does not respect Artin-Tate motives: we will see in \refsa{enoughpoints} that any Artin-Tate motive $M$ of the form $M = \eta_* M_\eta$, where $M_\eta$ is an  Artin-Tate motive over $F$, necessarily satisfies $M = 0$.

\pf
The stability of (Artin\nobreakdash-)Tate motives under $j^*$, $\eta^*$, $i_*$ and $i^*$, $f^*$ and---for Artin-Tate motives, under $f_*$---is immediate from the definition and base-change. For example, $i^* \phi_* \one (n) = \phi''_* \one(n) $. Here $\phi$ is any finite map over $S$ and $\phi''$ is its pullback along $i$. For the stability under $j_*$ it is sufficient to show $j_* \phi'_* \one$ is an Artin-Tate motive over $S$ for any finite flat map $\phi': V' \r S'$. Choose some finite flat (possibly non-regular) model $\phi : V \r S$ of $\phi'$, i.e., $V \x_S S' = V'$, so that $j^* \phi_* \one = \phi'_* \one$ is an Artin-Tate motive over $S'$. 
The localization triangle
$$j_* j^* \phi_* \one \r \phi_* \one \r i_* i^* \phi_* \one$$
and the above steps show that $j_* \phi_* \one$ is an Artin-Tate motive over $S$.

To see the stability under the Verdier dual functor $D$, it is enough to see that 
$$D(\phi_* \phi^* \one) = \phi_! \phi^! D(\one) = \phi_* \phi^! \one(1)[2].$$
is an Artin-Tate motive for any finite map $\phi : V \r S$ with reduced domain (\refbe{onepoint}). If $V$ is zero-dimensional, this follows from purity and the regularity of $S$. If not, there is an open (non-empty) immersion $j : S' \r S$ such that $V' := V \x_S S'$ is regular (for example, take $S'$ such that $V' / S'$ is \'etale). Let $i$ be the complement of $j$. We apply the localization triangle $i_* i^! \r \id \r j_* j^*$ to $\phi_* \phi^! \one$. By base-change we obtain
$$i_* \phi''_* \phi''^! i^! \one \r \phi_* \phi^! \one \r j_* \phi'_* \phi'^! j^* \one.$$
Here $\phi''$ and $\phi'$ is the pullback of $\phi$ along $i$ and $j$, respectively. By the regularity of $S$ and purity we have $i^! \one = \one(-1)[-2]$, so the left hand term is an Artin-Tate motive. The right one also is by purity. 
This shows the claim for $D$.

The stability under $f^!$, $i^!$, and $j_!$ now follow for duality reasons.


As for the stability under tensor products we note that $\phi_* \one \t \phi'_* \one = (\phi \x \phi')_* \one$ if $\phi$ and $\phi$ are (finite and) smooth, cf.\ \refeq{monoidal}. Using the localization triangle, it is easy to reduce the general case of merely finite maps $\phi$, $\phi'$ to this case. 
\xpf

\bem
\refle{respect} also holds for a similarly defined category of Artin-Tate motives over open subschemes $S$ of a smooth curve over a field.
\xbem

\lemm \mylabel{lemm_split}
Let $M \in \DATM(S)$ be any Artin-Tate motive. Then there is a finite map $f: V \r S$ such that $f^* M \in \DTM(S) \subset \DATM(S)$. We describe this by saying that $f$ \Def{splits} $M$.
\xlemm
\pf
As $f^*$ is triangulated, this statement is stable under triangles (with respect to $M$), and also under direct sums and summands. Therefore, we only have to check the generators, i.e., $M = \phi_* \one(n)$ with $\phi: S' \r S$ a finite map with reduced domain. The corresponding splitting statement for Artin-Tate motives over finite fields is well-known. Therefore, by localization, it is sufficient to find a splitting map $f$ after replacing $S$ by a suitable small open subscheme, so we may assume $\phi$ \'etale. We first assume that $\phi$ is moreover Galois of degree $d$, i.e., $S' \x_S S' \cong S'^{\sqcup d}$, a disjoint union of $d$ copies of $S'$. In that case one has $\phi^* \phi_* \one = \one^{\oplus d}$ by base-change, so the claim is clear. In general $\phi$ need not be Galois, so let $S''$ be the normalization of $S$ in some normal closure of the function field extension $k(S') / k(S)$. Both $\mu: S'' \r S$ and $\psi: S'' \r S'$ are generically Galois. By shrinking $S$ we may assume both are Galois. From $\Hom(\one_{S'}, \psi_* \one_{S''}) = \Hom(\one_{S''}, \one_{S''}) = \Q$ and $\Hom(\psi_* \one_{S''}, \one_{S'}) = \Hom(\one_{S''}, \psi^! \one_{S'}) = \Hom(\one_{S''}, \one_{S''}) = \Q$ we see that $\one_{S'}$ is a direct summand of $\psi_* \one_{S''}$. Therefore $\mu^* \phi_* \one_{S'}$ is a summand of $\mu^* \phi_* \psi_* \one_{S''} = \mu^* \mu_* \one_{S''} = \one^{\oplus \deg S'' / S}$, a Tate motive. 
\xpf

\section{The motivic $t$-structure}\mylabel{sect_motivict}
In this section, we establish the motivic $t$-structure on the category of Artin-Tate motives over $S$ (\refth{motivictstructure}). It is obtained by the standard gluing procedure, applied to the $t$-structures on Artin-Tate motives over finite fields and on a subcategory $\tilde T(S') \subset T(S')$ for open subschemes $S' \subset S$. Under the analogy of mixed (Artin-Tate) motives with perverse sheaves, the objects in the heart of the $t$-structure on $\tilde T(S')$ correspond to sheaves that are locally constant, i.e., have good reduction. We refer to \cite[Section 1.3.]{BBD} for generalities on $t$-structures.

\defi (compare \cite[Def. 1.1]{Levine:TateMotives}) \mylabel{defi_Tab}
For $-\infty \leq a \leq b \leq \infty$, let \Demm{Tab}{\tilde T_{[a, b]}} denote the smallest triangulated subcategory of $T(S)$ containing direct factors of $\phi_* \one(n)$, $a \leq -2n \leq b$, where $\phi: S' \r S$ is a finite \emph{\'etale} map. For Tate motives, $\phi$ is required to be the identity map. (We will not specify this restriction \emph{expressis verbis} in the sequel.) Furthermore, $\tilde T_{[a, a]}$ and $\tilde T_{[-\infty, \infty]}$ are denoted \Demm{Ta}{\tilde T_a} and \Demm{T}{\tilde T}. If it is necessary to specify the base, we write $\tilde T_{[a, b]}(S)$ etc.
\xdefi

We need the following vanishing properties of $K$-theory of number fields, related Dedekind rings and finite fields up to torsion. In order to weigh the material appropriately, it should be said that the content of the ``lemma'' below is the only non-formal part of the proofs in this paper, and all complexity occurring with Artin-Tate motives ultimately lies in these computations.

\lemm (Borel, Soul\'e, Quillen) \mylabel{lemm_morphisms}
Let $\phi: S' \r S$ and $\psi: V \r S$ be two finite maps with zero-dimensional domains.
\eqnarr
\Hom_S(\phi_* \one, \psi_* \one(n)[m]) & = & 
\left \{
\begin{array}{ll}

\text{finite-dimensional} & n=m=0 \\
0 & \text{else.} 
\end{array} \right. 
\xeqnarr

Let now $\phi: S' \r S$ and $\psi: V \r S$ be two finite \'etale maps over $S$. Then
%

\eqnarr
\Hom_S(\phi_* \one, \psi_* \one(n)[m]) & = & 
\left \{
\begin{array}{ll}
\text{finite-dimensional} & n = m = 0 \\
\text{finite-dimensional} & m = 1,\ n\text{ odd and positive} \\
0 & \text{else.} 
\end{array} \right. 
\xeqnarr

\xlemm
\pr
By \refeq{morphismsKTheory}
$$\Hom_V(\one, \one(q)[p]) \cong K_{2q-p}(V)_\Q^{(q)},$$ 
for a regular scheme $V$. For the first statement, we may assume that $S'$ and $V$ are finite fields. Then the statement follows from adjunction, base-change, purity and 
$$K_n(\Fq) = \left \{ 
\begin{array}{ll}
\mu_{q^i -1} & n = 2i-1, i > 0 \\
0 & n = 2i, i > 0 \\
\Z & n=0
\end{array} \right. $$
(Quillen \cite{Quillen:FiniteField}). 
$K$-theory of Dedekind rings $R$ whose quotient field is a number field is known (up to torsion) by Borel's work. The relation to $K$-theory of number fields is given by an exact sequence (due to Soul\'e \cite[Th. 3] {Soule:Anneaux}; up to two-torsion) for $n > 1$
$$0 \r K_n(R) \stackrel{\eta^*}{\lr} K_n(F) \r \oplus_{\pp} K_{n-1}(\Fpp) \r 0.$$
Here $\eta: \Spec F \r \Spec R$ is the generic point and the direct sum runs over all (finite) primes in $R$. Also, $K_0(R) = \Z \oplus Pic(R)$ and $K_1(R) = R^\x$. In particular, for all $n$ and $m$, $K_n(R)_\Q^{(m)}$ vanishes when $K_n(F)_\Q^{(m)}$ vanishes, since $\eta^*$ respects the Adams grading. One has the following list (see e.g.\ \cite{Weibel:KTheoryIntegers})
$$
K_{2q-p}(F)^{(q)}_\Q = \left \{ 
\begin{array}{ll}
0 & q < 0 \\
0 & q = 0, p \neq 0 \\
\Q & q = p = 0 \\
0^{BS} & q > 0, p \leq 0 \\
0 & q > 0, \text{even}, p=1 \\
F^\x \t_\Z \Q & q = p =1\\
\Q^{r_1 + r_2}  & q > 1, q \equiv 1 \mod{4}, p = 1 \\
\Q^{r_2}  & q > 0, q \equiv 3 \mod 4, p = 1 \\
0 & q > 0, p > 1
\end{array} \right. $$
As usual, $r_1$ and $r_2$ are the numbers of real and pairs of complex embeddings of $F$, respectively.
(The agreement of $K_{2q-1}(F)$ and $K_{2q-1}(F)^{(q)}$ for odd positive $q$ is not mentioned in \lc) The spot marked $0^{BS}$ is referred to as \Def{Beilinson-Soul\'e vanishing} (see e.g.\ \cite{Levine:TateMotives}). As first realized by Levine (\lc), this translates into the non-existence of morphisms in the ``wrong'' direction with respect to the motivic $t$-structure.

For the last claim, put $V' = V \x_S S'$:
$$\xymatrix
{
V' \ar[r]^{\phi'} \ar[d]^{\psi'} &
V \ar[d]^\psi \\
S' \ar[r]^\phi &
S.
}
$$
To save space, we omit the twist and the shift in writing the $\Hom$-groups. We have
$$\Hom_S(\phi_* \one, \psi_* \one) = \Hom_{S'}(\one, \phi^! \psi_* \one) = \Hom_{S'}(\one, \psi'_* \phi'^! \one) = \Hom_{V'}(\one, \phi'^! \one).$$
Now, $V'$ is (affine and) \'etale over $V$, so $\phi'^! \one = \phi'^* \one = \one$ and we are done in that case by the above vanishings of $K$-theory up to torsion. 
\xproof


The following lemma is a variant of \cite[Lemma 1.2]{Levine:TateMotives}, \cite[Lemma 1.9]{Wildeshaus:ATM} and can be proven by faithfully imitating the technique in \lcs

\lemm \mylabel{lemm_auxiliaryt}
For any $-\infty \leq a < b \leq c \leq \infty$, $(\tilde T_{[a, b-1]}, \tilde T_{[b, c]})$ is a $t$-structure on $\tilde T_{[a, c]}$. 
\xlemm

\defi
The resulting truncation and cohomology functors are denoted $F_{\leq b}$ and $F_{>b}$ and $\gr_b^F$, respectively.  
\xdefi

The following definition is modeled on \cite[Def. 1.4]{Levine:TateMotives}. We also refer to \cite[Section 2.1.3]{Ayoub:Six1} for a general way (due to Morel) of constructing a $t$-structure starting from a given set of generators. For any odd integer $n$ set $\one (n/2) := 0$, for notational convenience.

\defi \mylabel{defi_tstructure} 
Let $S$ be an open subscheme of $\SpecOF$. 
Let $\tilde T^{\geq 0}_a(S)$ ($\tilde T^{\leq 0}_a(S)$) be the full subcategory of $\tilde T_a(S)$ (\refde{Tab}) generated by 
$$\phi_* \one \left( -\frac a 2 \right)[n+1]$$
for any $n \leq 0$ ($n \geq 0$, respectively), and any finite \'etale map $\phi$. 
``Generated'' means the smallest subcategory containing the given generators stable under isomorphism, finite direct sums, summands and $\cone(\phi)[-1]$ ($\cone(\phi)$, resp.) for any morphism $\phi$ in $\tilde T_a^{\geq 0}(S)$ ($\tilde T_a^{\leq 0}(S)$, respectively).

For any $-\infty \leq a \leq b \leq \infty$, let $\tilde T_{[a,b]}^{\geq 0}(S)$ be the triangulated subcategory generated by objects $X$, such that for all $a \leq c \leq b$, $\gr_c^F(X) \in \tilde T_{c}^{\geq 0}(S)$ and similarly for $\tilde T_{[a,b]}^{\leq 0}(S)$. For $a = -\infty$ and $b= \infty$ we simply write $\tilde T^{\leq 0}(S)$, $\tilde T^{\geq 0}(S)$. We may omit $S$ in the notation, if no confusion arises.
\xdefi

In particular $\one(-a/2)[1] \in \tilde T^0_a(S)$. This shift is as in the situation of perverse sheaves \cite{BBD}, \cite[\fcohoSectionThree]{Scholbach:fcoho}. Before stating and proving the existence of the motivic $t$-structure, we need some preparatory steps. Levine has established the existence of the motivic $t$-structure on Tate motives over number fields and finite fields \cite[Theorem 1.4.]{Levine:TateMotives}. 
This has been generalized to Artin-Tate motives by Wildeshaus \cite[Theorem 3.1]{Wildeshaus:ATM}. We briefly recall these precursor statements. Let $K$ be either a finite field or a number field. 
For any $-\infty \leq a \leq b \leq \infty$, let $T_{[a,b]}(K)$ be the triangulated subcategory of $T(K)$ generated by $\one(n)$ with $a \leq -2n \leq b$ (Tate motives) and direct summands of $\phi_* \one(n)$, $\phi : \Spec K' \r \Spec K$ a finite map (Artin-Tate motives, respectively). For any $a \leq c < b$, the datum $\left ( T_{[a,c]}, T_{[c+1, b]} \right )$ forms a $t$-structure on $ T_{[a,b]}$. Let $\gr^F_*$ be the cohomology functor corresponding to that $t$-structure. Write $T_a(K)$ for $T_{[a,a]}(K)$ and let $T^{\geq 0}_a(K)$ and $T^{\leq 0}_a(K)$ be the subcategories of $T_a(K)$ generated by $\one(-a/2)[n]$ with $n \leq 0$ and $n \geq 0$, respectively. Here, ``generated'' has the same meaning as  in \refde{tstructure}. Let $T^{\geq 0}_{[a,b]}$ and $T^{\leq 0}_{[a,b]}$ be the subcategories of $T_{[a,b]}$ of objects $X$ such that all $\gr^F_c X \in T^{\geq 0}_c$ ($\gr^F_c X \in T^{\leq 0}_c$, respectively) for all $a \leq c \leq b$. Then, $\left ( T^{\leq 0}_{[a,b]}(K), T^{\geq 0}_{[a,b]}(K) \right)$ is a non-degenerate $t$-structure on $T_{[a,b]}$. 

The following well-known fact is a consequence of vanishing of all $K$-theory groups of finite fields except for $K_0(\Fpp)^{(0)}_\Q$, see \refle{morphisms}.

\lemm \mylabel{lemm_tstructureFp}
Let $\pp$ be a closed point in $S$ with residue field $\Fpp$. The inclusions $T_a(\Fpp) \subset T(\Fpp)$ induce an equivalence of categories
$$\bigoplus_{a \in \Z} T_a(\Fpp) = T(\Fpp).$$
\xlemm

There are canonical equivalences of categories 
$$T(Z) := 
\bigoplus_{\pp \in Z, a \in \Z} T_a (\Fpp) = 
\bigoplus_{\pp, a} 
\D^\bound(\underline \Q[\text{Perm}, \Gal(\Fpp)]) 
= \bigoplus_{\pp, a} \underline \Q[\text{Perm}, \Gal(\Fpp)]^{\Z-\text{graded}}
$$
Here and in the sequel $\underline \Q[\text{Perm}, \Gal(\Fpp)]$ denotes finite-dimensional rational permutation representations of the absolute Galois group. By means of that equivalence, $T(Z)$ is endowed with the obvious $t$-structure. 
The heart $T^0_a (\Fpp) = T^{\leq 0}_a (\Fpp) \cap T^{\geq 0}_a (\Fpp)$ is semisimple and consists of direct sums of summands of $\phi_* \one(a)$, $\phi$ finite. 

We now provide the motivic $t$-structure on $\tilde T(S)$, which stems from the one on $T(F)$. The two together will then be glued to give the $t$-structure on $T(S)$. Recognizably, the following is again an adaptation of Levine's proof of the $t$-structure on Tate motives over number fields. 

\satz \mylabel{satz_prep1}
For any $-\infty \leq a \leq b \leq \infty$, $\left( \tilde T_{[a,b]}^{\leq 0}, \tilde T_{[a,b]}^{\geq 0} \right)$ is a non-degenerate $t$-structure on $\tilde T_{[a,b]}(S)$ (Definitions \ref{defi_Tab}, \ref{defi_tstructure}). The functor $\eta^*[-1] : \tilde T_{[a,b]}(S) \r T_{[a, b]}(F)$ is $t$-exact.

Any motive in $\tilde T^0_a(S)$ is a finite direct sum of summands of motives $\phi_* \one(-a/2)[1]$ with $\phi$ finite \'etale. The closure of the direct sum of the $\tilde T^0_a(S)$, $a \in \Z$, under extensions (in the abelian category $\tilde T^0(S)$) is $\tilde T^0(S)$.
\xsatz

\pf
We may assume that $a$ and $b$ are finite, since 
$$\tilde T(S) = \bigcup_{-\infty < a \leq b < \infty} \tilde T_{[a, b]}(S)$$ 
and the inclusion functors given by the identity between the various $T_{[-, -]}$ are exact. 

The proof proceeds by induction on $b-a$. The case $b=a$ is treated as follows: the category $\tilde T_a := \tilde T_a(S)$ is generated by $\phi_* \one(-a/2)[n]$, $n \in \Z$, $\phi$ \'etale and finite. The functor $\eta^*[-1](a/2): \tilde T_a(S) \r T_0(F)$ is fully faithful. To see this it suffices to remark $\Hom_S(\phi_* \one(-a/2)[n+1], \psi_* \one(-a/2)[n'+1]) = \Hom_F({\phi_\eta}_* \one[n], {\psi_\eta}_* \one[n'])$, for any finite \'etale maps $\phi$ and $\psi$ with generic fiber $\phi_\eta$ and $\psi_\eta$. This equality follows from the $K$-theory computations, see the proof of \refle{morphisms}. Therefore, the image of $\eta^*[-1](a/2)$ is a triangulated subcategory of $T_0(F)$ which contains the generators of $T_0(F)$, so the functor establishes an equivalence between $\tilde T_a(S)$ with the derived category of finite-dimensional rational permutation representations of $\Gal(F)$ by \cite[3.4.1]{Voevodsky:TCM}. Hence $\tilde T_a(S)$ carries a non-degenerate $t$-structure.

The remainder of the proof is done as in Levine's proof. One shows 
\begin{equation}
\Hom \left(\tilde T_{[a+1,b]}^{\leq 0}, \tilde T_{c}^{\geq 0} \right) = 0 \label{eqn_Hom}
\end{equation}
for any $c \leq a$. This reduces to the Beilinson-Soul\'e vanishing. Then the $t$-structure axioms follow for formal reasons. 

The exactness of $\eta^*[-1]$ is obvious from the definitions. The statement about the heart $\tilde T_a^0$ is done as follows: the exact functor $\eta^*[-1](a/2)$ identifies $\tilde T^0_a(S) = \tilde T^{\geq 0}_a(S) \cap \tilde T^{\leq 0}_a(S)$ with the semi-simple category $T^0_0(F) = \underline \Q[\text{Perm}, \Gal(F)]$. We claim that for any object $X \in \tilde T_a(S)$, all $\pH^n(X)$ are direct summands of sums of motives $\phi_* \one(-a/2)[1]$, $\phi$ finite and \'etale. This claim does hold for the generators of $\tilde T_a(S)$. We now show that the condition is stable under triangles, which accomplishes the proof of the claim and thus the proof of the statement. Let $A \r X \r B$ be a triangle in $\tilde T_a(S)$ such that $A$ and $B$ satisfy the claim. The long exact cohomology sequence 
$$\dots \r \pH^{n-1} B \stackrel{\delta^{n-1}}{\lr} \pH^n A \r \pH^n X \r \pH^n B \stackrel{\delta^{n}}{\lr} \pH^{n+1} A \r \dots$$
yields the short exact sequence in $\tilde T^0_a(S)$
$$0 \r \coker \delta^{n-1} \r \pH^n X  \r \ker \delta^n \r 0.$$
By the semi-simplicity of $\tilde T^0_a(S)$ (this is the key point!), the sequence splits and there is a non-canonical isomorphism $\pH^n X \cong \coker \delta^{n-1} \oplus \ker \delta^n$ and $\coker \delta^{n-1}$ and $\ker \delta^n$ are direct summands of $\pH^n A$ and $\pH^n B$, respectively.

For the statement concerning $\tilde T^0(S)$ one uses the finite exhaustive $F$-filtration of any $X \in \tilde T^0(S)$:
$$0 = F_a X \subset F_{[a, a+1]} X \subset \dots \subset F_{[a, b]} X = X.$$
The successive quotients $\gr_*^F X$ of that chain are in $\tilde T^0_*(S)$, since truncations with respect to the $t$-structure related to $F$ are exact with respect to the motivic $t$-structure, by definition. Thus the claim about $\tilde T^0(S)$ follows.
\xpf

\theo \mylabel{theo_motivictstructure}
The motivic $t$-structures on $T(Z)$ and $\tilde T(S')$ glue to a non-degenerate $t$-structure on the category $T(S)$ of (Artin\nobreakdash-)Tate motives over $S$ (\refde{DATM}). It is called \Def{motivic $t$-structure}. Here $S'$ runs through open subschemes of $S$ and $Z := S \backslash S'$. 
\xtheo
\pr
We apply the gluing procedure of $t$-structures of \cite[Theorem 1.4.10]{BBD}: for any open subscheme $j: S' \subset S$, we write $T_{S'}(S)$ for the full triangulated subcategory of objects $X \in T(S)$ such that $j^* X \in \tilde T(S') \subset T(S')$. Let $i : Z' \r S$ be the closed complement of $j$. 
Put
$$T^{\leq 0}_{S'}(S) := \{ X \in T_{S'}(S), j^* X \in \tilde T^{\leq 0}(S'), i^* X \in T^{\leq 0}(Z') \},$$
$$T^{\geq 0}_{S'}(S) := \{ X \in T_{S'}(S), j^* X \in \tilde T^{\geq 0}(S'), i^! X \in T^{\geq 0}(Z') \}.$$
The assumptions of the gluing theorem, \cite[1.4.3]{BBD}, namely the existence of $i_*$, $i^*$, $i^!$, $j_*$, $j_!$, $j^*$ satisfying the usual adjointness properties, $j^* i_* = 0$, localization sequences and full faithfulness of $i_*$, $j_!$ and $j_*$ are met, since they are in the surrounding categories of geometric motives, cf. \refsect{preliminaries}, and the stability results of \refsect{definitionstability}. Thus, the above defines a $t$-structure on $T_{S'}(S)$.

The field $F$ is of characteristic zero, so any finite map $\phi: V \r S$ with $V$ reduced and one-dimensional is generically \'etale. This implies $T(S) = \cup_{S' \subset S} T_{S'}(S)$. We set 
$$T^{\geq 0}(S) := \bigcup_{S' \subset S} T^{\geq 0}_{S'}(S)$$
and dually for $T^{\leq 0}(S)$. The $t$-structure axioms on $T(S)$ and the non-degeneracy are implied by the exactness of the identical inclusion $T_{S'}(S) \r T_{S''}(S)$ for any $S'' \subset S'$. 

To see the exactness of the identity, let $j'' : S'' \subset S$ and $i'': Z'' \subset S$ be its complement. Let $X \in T^{\leq 0}_{S'}(S)$. It is clear that $j''^* X \in \tilde T^{\leq 0}(S'')$. Let us check $i''^* X \in T^{\leq 0}(Z'')$. The pullback $i''^* X $ decomposes as a direct sum parametrized by the points of $Z''$ and we only have to deal with the points that are not contained in $Z'$. Let $p: \SpecFpp \r S$ be such a point; it factors over $S'$: $p = j \circ q$, where $q: \SpecFpp \r S'$ is the same point as $p$. Thus $p^* X = q^* j^* X \in q^* \tilde T^{\leq 0}(S')$. The containment $ q^* \tilde T^{\leq 0}(S') \subset T^{\leq 0}(\SpecFpp)$ follows from $q^* \tilde T_a^{\leq 0}(S') \subset T_a^{\leq 0}(\SpecFpp)$, since $q^*$ clearly commutes with the $F$-truncation functors belonging to the auxiliary $t$-structure. To see the latter containment, it suffices to check the generators (in the sense of \refde{tstructure}) of $\tilde T_a^{\leq 0}(S')$, that is, it is sufficient to remark
$$q^* \phi_* \one(-a/2)[n+1] = \phi'_* \one (-a/2)[n+1] \in T_a^{\leq -1}(\SpecFpp) \subset T_a^{\leq 0}(\SpecFpp),$$ 
where $n \geq 0$ and $\phi$ is a finite \'etale map with pullback $\phi'$. 
This shows that the identity is left-exact. The right-exactness is done dually. 
\xpf

\lemm \mylabel{lemm_summands} 
The category $\DTM(S)$ agrees with the triangulated category generated by $\one(n)$, $i_* \one(n)$. 
\xlemm
\pf
Let $M \in \DTM(S)$. Pick an open subscheme $j: S' \subset S$ with complement $i: Z \subset S$ such that $j^* M \in \tilde T(S')$. Any object in $T(Z)$ is isomorphic to a direct sum of motives $\one_{\Fpp}(a)[b]$, $\pp \in Z$, since $\one_\Fpp$ does not have proper direct summands. Any object in $\tilde T^0_{-2a}(S')$ is a direct sum of motives $\one(a)[1]$ for the same reason. Any object in $\tilde T^0(S')$ is obtained by taking repeated extensions starting with such objects. Thus $\tilde T(S')$ is the triangulated category generated by $\one(a)$, $a \in \Z$. The localization triangle $i_* i^! M \r M \r j_* j^* M$ settles the lemma.
\xpf

\section{Mixed Artin-Tate motives} \mylabel{sect_exactness}

\defi \mylabel{defi_mixedAT}
The heart $T^0(S)$ of the motivic $t$-structure is called the category of \emph{mixed (Artin\nobreakdash-)Tate motives} over $S$, denoted $\MTM(S)$ and $\MATM(S)$, respectively. The cohomology functors belonging to the motivic $t$-structure are denoted $\pH^*$.
\xdefi

We now study the categories of mixed Tate motives over $S$ in some detail. The key is \refth{exactness} below, establishing exactness properties of pullback and pushforward functors along closed and open immersions. The exactness axioms for mixed motives over number rings (see \cite[\fcohoSectionFour]{Scholbach:fcoho}) are modeled on this theorem. Of course, the theorem is an Artin-Tate motivic analog of a similar fact about perverse sheaves \cite[Prop. 1.4.16, 4.2.4.]{BBD}, suggesting that the theory of perverse sheaves is to some extent quite formal. \refsa{cohomdimAT} calculates the cohomological dimension of mixed (Artin\nobreakdash-)Tate motives. We obtain an equivalence $\DTM(S) \cong \D^\bound(\MTM(S))$, using a result of Wildeshaus, and likewise for Artin-Tate motives. Finally, we do a first step into (Artin-Tate) motivic sheaves, in \refsa{enoughpoints}.

All exactness statements below are with respect to the motivic $t$-structure of \refth{motivictstructure}. Recall from \refle{respect} that the functors discussed below do preserve (Artin\nobreakdash-)Tate motives. For brevity, we write $T^{[a, b]}$ for the full subcategory of objects $M$ satisfying $\pH^n M=0$ for all $n < a$ and $n > b$. We say that a triangulated functor $F$ between categories of Artin-Tate motives has \Def{cohomological amplitude} $[a, b]$ if $F(T^0)$ is contained in $T^{[a, b]}$. Note that $F$ is right exact iff $b \leq 0$ and left exact iff $a \geq 0$. 

\theo \mylabel{theo_exactness}
Let $j : S' \r S$ be an open immersion, $i: Z \r S$ a closed immersion with $\dim Z = 0$. Finally, let $f: V \r S$ be a finite map with regular one-dimensional domain.
\begin{enumerate}[(i)]
\item \mylabel{item_Verdier}
The Verdier duality functor $D$ is exact in the sense that it maps $T^{\geq 0}$ to $T^{\leq 0}$ and vice versa. Therefore, it induces an endofunctor on $T^0(S)$. 
\item \mylabel{item_ij}
The functors $j_*$, $j_!$, $j^*$, as well as $i_* = i_!$ are exact. 
\item \mylabel{item_iupperstar}
The functor $i^*$ has cohomological amplitude $[-1, 0]$. Dually, $i^!$ has cohomological amplitude $[0, 1]$.
\item \mylabel{item_f}
The functor $f_* = f_!$ is exact. The cohomological amplitude of $f^*$ and $f^!$ is $[-1, 0]$ and $[0, 1]$, respectively. If $f$ is also \'etale, $f^* = f^!$ is exact.  
\item \mylabel{item_eta}
The functor $\eta^*[-1]: T(S) \r T(\Spec F)$ is exact. 
\end{enumerate}
\xtheo
\pr
\refit{Verdier}: This is clear from the definitions of the $t$-structures on $T(S)$, $\tilde T(S')$ and $T(Z)$, for open and closed subschemes $S'$ and $Z$ of $S$, respectively. Notice that this requires putting $\one[1]$ in degree $0$. 

\refit{ij}: The following exactness properties are immediate from the definition: $j^*$ and $i_*$ are exact, $j_*$ and $i^!$ are left-exact and $j_!$ and $i^*$ are right-exact. For example, let us show the left-exactness of $j_*$. Given some motive $M \in T^{\geq 0}(S')$, we have to show $j_* M \in T^{\geq 0}(S)$. Let $j_1: S_1 \subset S'$ be an open immersion such that $j_1^* M \in \tilde T^{\geq 0}(S_1)$. Let $i_1$ be the immersion of $Z_1 := S' \backslash S_1$ into $S'$, then $i_1^! M \in T^{\geq 0}(Z_1)$. The situation is as follows:
$$\xymatrix{
& & Z_1 \ar[dl]^{i_1} \ar[dr] \\
S_1 \ar[drr] \ar[r]^{j_1}&
S' \ar[dr]^j &
&
S \backslash S_1 \ar[dl]^i &
\\
& & S
}
$$
Now $(j \circ j_1)^* j_* M = j_1^* M \in T^{\geq 0}(S_1)$. Let $i: S \backslash S_1 \r S$ be the complement of $j \circ j_1$. Then $i^! j_*M$ is supported only in $Z_1$, where it agrees with $i_1^! M$. This shows $j_* M \in T^{\geq 0}(S)$.

To prove \refit{iupperstar} we first show 
\eqn \mylabel{eqn_claimexactness}
i^* j_* \tilde T^0(S') \subset T^{[-1,0]}(Z)
\xeqn 
for any two complementary immersions $i: Z \r S$ (closed) and $j: S' \r S$ (open). By \refsa{prep1}, $\tilde T^0(S)$ is generated by means of direct sums, summands and extensions by $\phi_* \one(n)[1]$, where $n \in \Z$ is arbitrary and $\phi$ is finite and \'etale. For any short exact sequence
$$0 \r A \r X \r B \r 0$$
in $\tilde T^0(S)$, such that $i^* j_* A \in T^{[-1,0]}(Z)$ and $i^* j_* B \in T^{[-1,0]}(Z)$, it follows $i^* j_* X \in T^{[-1,0]}(Z)$. This uses the non-degeneracy of the motivic $t$-structure on $Z$. A similar remark applies to direct summands and sums. Therefore we only have to check that the generators $X$ of $\tilde T^0(S')$ are mapped to $T^{[-1,0]}(Z)$ under $i^* j_*$. Thus, let $X = \phi_* \one(n)[1]$. We have a localization triangle in $T(Z)$
$$i^* \phi_* \one(n)[1] \r i^* j_* j^* \phi_* \one(n)[1] = i^* j_* \phi'_* \one(n)[1] \r i^! \phi_* \one(n)[2] \r i^* \phi_* \one(n)[2].$$
Here $\phi'$ is the pullback of $\phi$ along $j$.
The first term is in degree $-1$. The third term is in degree $0$ by absolute purity (see \refsect{preliminaries}), using the regularity of $S$. The claim \refeq{claimexactness} is shown. 

We now show $i^* T^0(S) \subset T^{[-1,0]}(Z)$. Any $X \in T^0(S)$ is in some $T^0_{S'}(S)$ for sufficiently small $S'$. We shrink $S'$ if necessary to ensure that $S' \cap Z = \emptyset$. Let $j: S' \r S$ be the open immersion and let $p: W \r S$ be its closed complement. There is a triangle
$$p^! X \r p^* X \r p^* j_* j^* X \r p^! X[1].$$
We know that $p^!$ ($p^*$) is left-exact (right-exact), that is to say, the first (second) term is in degrees $\geq 0$ ($\leq 0$, respectively). By assumption $j^* X \in \tilde T^0(S')$, so $p^* j_* j^* X \in T^{[-1,0]}(W)$, as was shown above. As the $t$-structure on $W$ is non-degenerate $p^* X$ is in degrees $[-1, 0]$. As $W$ is the disjoint union of $Z$ and some more (finitely many) closed points, this also shows $i^* X \in T^{[-1,0]}(Z)$.

Let now $i: Z \r S$ and $j: S' \r S$ be complementary. We claim $i^* j_* T^0(S') \subset T^{[-1,0]}(Z)$. Given an object $X \in T^0(S')$, there is some open immersion $j': S'' \r S'$ such that $j'^* X \in \tilde T^0(S'')$. We have $i^* j_* X = i^* j_* j'_* j'^* X$. The motive $i_* i^* j_* j'_* j'^* X $ is a direct summand of $p_* p^* (j \circ j')_* j'^*X$, where $p$ is the complement of $j \circ j'$. By the above, $p^* (j \circ j')_* j'^*X \in T^{[-1,0]}(Z)$, so the full faithfulness and exactness of $p_*$ implies the claim. \refit{iupperstar} is shown.

The cohomological amplitude of $i^* j_*$ implies the exactness of $j_*$: given a mixed (Artin\nobreakdash-)Tate motive $M \in T^0(S')$, the terms in the localization triangle 
$$j_! M \r j_* M \r i_* i^* j_* M$$
are in degrees $\leq 0$, $\geq 0$ and $[-1, 0]$, respectively, by the above. From the non-degeneracy of the $t$-structure we see that $j_* M$ is then in degree $0$. This implies the exactness of $j_*$ by the non-degeneracy of the $t$-structure. 

The exactness of $j_!$ follows by Verdier duality, as does the cohomological amplitude of $i^!$. Thus, \refit{ij} is shown.

\refit{f}: It is easy to see that $f^*: \tilde T(S) \r \tilde T(V)$ is exact. Using this and the localization triangles, one sees that $f^*$ has cohomological amplitude $[-1, 0]$ and dually for $f^!$. By a general criterion on $t$-exactness of adjoint functors \cite[1.3.17]{BBD}, the adjunctions $f^* \leftrightarrows f_* = f_! \leftrightarrows f^!$ imply that $f_*$ is exact. If $f$ is \'etale then $f^! = f^*$, so that their exactness is clear in that case, too.  

\refit{eta}: This follows from the exactness of $j^* : T(S) \r T(S')$ and the exactness of $\eta'^*[-1]: \tilde T(S') \r T(\SpecF)$ (\refsa{prep1}), where $\eta'$ is the generic point of $S'$.
\xpf

\defi (\cite[1.4.22]{BBD}, see \cite[\fcohoSectionFour]{Scholbach:fcoho} for the motivic case)
Let $j: S' \r S$ be an open immersion. For any mixed (Artin\nobreakdash-)Tate motive $M$ over $S'$, put $j_{!*} M = \im j_! M \r j_* M$. This is called the \Def{intermediate extension} of $M$ along $j$. 
\xdefi

The image is taken in the (abelian) category of mixed (Artin\nobreakdash-)Tate motives over $S$, using the exactness of $j_!$ and $j_*$. Thereby, $j_{!*}$ is a (non-exact) functor $T^0(S') \r T^0(S)$. Given any mixed motive $M$ over $S$, such that $i^! M$ is concentrated in cohomological degree $-1$ (as opposed to the general range $[-1, 0]$), and such that $i^* M$ is in degree $+1$, there is a canonical isomorphism
\eqn
\mylabel{eqn_intermediate}j_{!*} j^* M = M.
\xeqn
In particular, this applies to $M \in \tilde T^0(S)$, such as $M = \one[1]$. Moreover, taking the intermediate extension commutes with compositions of open immersions. These features will be used below, see \lcs for a proof. The reader may want to check that that proof only uses the motivic $t$-structure and exactness properties of $i^!$ etc., which are established by Theorems \ref{theo_motivictstructure}, \ref{theo_exactness}.

\satz \mylabel{satz_cohomdimAT}
The cohomological dimension of $\DTM(S)$ and $\DATM(S)$ is one and two, respectively. 
\xsatz

\pf
We have to show $\Hom(M, M'[n]) = 0$ for any mixed motives $M$, $M'$ over $S$ and $n > 1$ (Tate) and $n > 2$ (Artin-Tate). Let $j: S' \r S$ be an open immersion such that $j^* M$, $j^* M' \in \tilde T^0(S')$. Let $i$ be the complementary closed immersion of $j$. In the sequel we write $(-,-)^n$ for $\Hom(-,-[n])$ for brevity. 

The case $n \geq 3$ is done as follows: the localization triangle \refeq{localization2} for $M'$ and adjunction gives a long exact sequence
$$(\underbrace{i^* M}_{[-1,0]}, \underbrace{i^! M'[n]}_{[-n, -n+1]})^0 \r (M, M')^n \r (M, j_* j^* M')^n \r (\underbrace{i^* M}_{[-1,0]}, \underbrace{i^! M'[n+1]}_{[-n-1, -n]})^0$$
We have written the cohomological degrees of the motives underneath, using the cohomological range of $i^*$ and $i^!$. The cohomological dimension zero of (Artin\nobreakdash-)Tate motives over finite fields makes the outer terms vanish. Similar vanishings will be used below without further discussion. Hence we only have to look at $(j^* M, j^* M')^n$, i.e., we may assume $M$ and $M' \in \tilde T^0(S)$. In that case one reduces (exactly as below) to $M = \phi_* \one(a)[1]$ and $M = \phi'_* \one(a')[1]$, where $\phi$ and $\phi'$ are finite and \'etale. 
In that case the vanishing is given by \refle{morphisms}.

The vanishing in the case $n=2$ for Tate motives needs a more involved localization argument. A similar reasoning for Artin-Tate motives fails---the difference is because the motives $\one(n)[1]$, which generate $\tilde T^0(S)$ in the case of Tate motives, have good reduction at all places by absolute purity.

The localization triangle for $M'$ gives an exact sequence
$$(M, j_! j^* M')^2 \r (M, M')^2 \r (M, i_* i^* M')^2 = (\underbrace{i^* M}_{[-1,0]}, \underbrace {i^* M'[2]}_{[-3,-2]})^0 = 0.$$
Therefore, in order to show that the middle term vanishes, we may replace $M'$ by $j_! j^* M'$. Similarly, we may replace $M$ by $j_* j^* M$. In particular $M \in j_* \tilde T^0(S')$, $M' \in j_! \tilde T^0(S')$. By \refsa{prep1}, $\tilde T^0(S')$ is generated by means of extension and direct summands by $\one(a)[1]$ where $a \in \Z$. The claim is stable under extensions and direct summands and sums  so that we may assume $M = j_* A$, $A := \one(a)[1]$, $M' = j_! A'$, $A' := \one(a')[1]$. Let $\tilde A := \one (a)[1] \in \tilde T^0(S)$ and define $\tilde A'$ similarly. We have $j^* \tilde A = A$ and similarly with $A'$. 

The triangle $j_* A' \r i_* i^* j_* A' \r j_! A'[1]$ maps to $j_* A' \r i_* \pH^0 i^* j_* A' \r (j_{!*} A')[1] = \tilde A[1]$. We apply $(\tilde A, -)^1$ to this map, which gives the last two exact rows in the diagram. The first exact row maps to the second via the adjunction map $\tilde A = j_{!*} A \r j_* A$.
$$\xymatrix{
(j_* A, j_* A')^1 \ar[r] \ar@2{-}[d] &
(j_* A, i_* i^* j_* A')^1 \ar[r] \ar@2{-}[d] &
(j_* A, j_! A')^2 \ar[d] \ar[r] &
0 \\
(\tilde A, j_* A')^1 \ar@2{-}[d] \ar[r] &
(\tilde A, i_* i^* j_* A')^1 \ar@2{-}[d] \ar[r] &
(\tilde A, j_! A)^2 \ar[d] \ar[r] &
0
\\
(\tilde A, j_* A')^1 \ar[r] &
(\tilde A, i_* \pH^0 i^* j_* A')^1  \ar[r] &
(\tilde A, \tilde A)^2 \ar[r] &
0
}
$$
The $=$ signs in the leftmost column are by adjunction and $j^* j_* A = j^* \tilde A = A$. The $=$ signs in the second column all use the adjunction $i^* \leftrightarrows i_*$ as well as the comological dimension zero of Tate motives over finite fields and cohomological amplitude of $i^*$, which imply 
$$(\underbrace{i^* j_* A}_{[-1, 0]}, \underbrace{i^* j_* A'[1]}_{[-2, -1]})^0 = (\pH^{-1} i^* j_* A, \pH^0 i^* j_* A')^0.$$
Applying $i^*$ to the triangle $i_* \pH^{-1} i^* j_* A \r j_! A \r j_{!*} A$ and using $i^* j_! = 0$ we see $(\pH^{-1} i^* j_* A, \pH^0 i^* j_* A')^0 = (i^* j_{!*} A, \pH^0 i^* j_* A')^1$. This justifies the upper $=$ in the second column. The lower $=$ in that column follows by the same argument. However, $(\tilde A, \tilde A')^2 = 0$, by vanishing of $K$-theory in the relevant range (see \refle{morphisms}).
\xpf


\theo \mylabel{theo_Wildeshaus}
For both Tate and Artin-Tate motives, the inclusion $T^0(S) \subset T(S)$ extends to a triangulated functor
\eqn 
\mylabel{eqn_DbTate}
\D^\bound(T^0(S)) \r T(S).
\xeqn
This functor is an equivalence of categories. 
\xtheo
\pf
The category $\DMgm(S)$ and thus the subcategories of (Artin\nobreakdash-)Tate motives embedd into some unbounded derived category $\D(\mathcal A)$, where $\mathcal A$ is an exact category. This implies the first statement by a general fact in homological algebra \cite[Theorem 1.1.]{Wildeshaus:fcat}. 
Indeed, the interpretation of $\DMgm(S)$ in terms of $h$-sheaves shows that (using the notation of \cite{CisinskiDeglise:Triangulated} and abbreviating $\Shvv$ for the category of $\Q$-linear sheaves with respect to the $h$-topology on the big site of schemes of finite type over $S$)
$$\DMgm(S) \cong 
\D_{\A}(\Shvv) \subset 
\D^\eff_{\A}(\category{Sp}(\Shvv)) \subset 
\D(\category{Sp}(\Shvv)).$$
More precisely, $\DMgm(S)$ identifies with the subcategory of $W_\Omega$-local objects in the middle category, which identifies with the subcategory of $W_{\A}$-local objects in the right hand category \cite[Sections 5.2, 5.3]{CisinskiDeglise:Triangulated}. 

The $t$-structure on $T(S)$ is bounded and non-degenerate, so it remains to show the full faithfulness of \refeq{DbTate} or equivalently that the map 
$$f_n: \Ext^n_{T^0}(M, M') \r \Hom_T (M, M'[n])$$
is an isomorphism for any $M$, $M' \in T^0(S)$. The general theory shows that $f_0$ and $f_1$ are isomorphisms and that $f_2$ is injective for all $M$ and $M'$. For Tate motives, $f_2$ is therefore an isomorphism, since the right hand side is zero by \refsa{cohomdimAT}. We next show that $f_2$ is an isomorphism for Artin-Tate motives. The motives $M$ and $M'$ are fixed, so there is some open embedding $j: S' \r S$ such that $j^* M$ and $j^* M'$ are in $\tilde T^0(S')$. Let $i$ be the complement of $j$. Consider the exact localization sequences
\eqn
\mylabel{eqn_exactseqproof}
0 \r i_* \pH^{-1} i^* M \stackrel{a}\r j_! j^* M \r K := \coker a \r 0
\xeqn
\eqn
\mylabel{eqn_exactseq2proof}
0 \r K \r M \r i_* \pH^{0} i^* M \r 0. 
\xeqn
We write ${}^n(-,-)$ for $\Ext^n$ and ${}_n(-,-)$ for $\Hom_T(-,-[n])$. \refeq{exactseqproof} induces a commutative diagram with exact rows
$$
\xymatrix{
{}^1(i_* \pH^{-1} i^* M, M') \ar[r] \ar@{=}[d] & 
{}^2(K, M') \ar[r] \ar@{>->}[d] &
{}^2(j_! j^* M, M') \ar@{>->}[d] \\
{}_1(i_* \pH^{-1} i^* M, M') \ar[r] &
{}_2(K, M') \ar[r] &
{}_2(j_! j^* M, M') = {}_2(j^* M, j^* M').
}
$$
The rightmost lower term is zero by the vanishings of $K$-theory (cf. the argument in the proof of \refsa{cohomdimAT}), so all vertical maps are isomorphisms. This and \refeq{exactseq2proof} yields a similar diagram:
$$
\xymatrix{
{}^2(i_* \pH^0 i^* M, M') \ar[r] \ar@{>->}[d] & 
{}^2(M, M') \ar[r] \ar@{>->}[d]^r &
{}^2(K, M') \ar@{=}[d] \ar[r] &
{}^3(i_* \pH^0 i^* M, M') \ar[d]
\\
{}_2(i_* \pH^0 i^* M, M') \ar[r] &
{}_2(M, M') \ar[r] &
{}_2(K, M') \ar[r] &
{}_3(i_* \pH^0 i^* M, M')
}
$$
The outer terms in the lower row vanish because the cohomological dimension of Artin-Tate motives over $\Fpp$ is zero and $i^!$ has cohomological amplitude $[0, 1]$. We now show that the rightmost upper term is zero. Altogether, this implies that $r$ is also surjective. We write $A := \pH^0 i^* M$; it is a mixed motive over $\Fpp$. Any element of the Yoneda-$\Ext$-group in question is represented by an exact sequence 
$$0 \r i_* A \r X_1 \stackrel{s}\r X_2 \r X_3 \r M' \r 0$$
in $\MATM(S)$. This extension is the image under the concatenation mapping 
$${}^2 (i_* A, \coker s) \x {}^1 (\coker s, M') \r {}^3 (i_* A, M')$$
The left hand factor is a subgroup of ${}_2 (i_* A, \coker s) = {}_2 (A, i^! \coker s) = 0$ (see above). Therefore, the extension above splits and we have shown that second $\Ext$-groups and $\Hom$-groups agree. 

This shows that the $\Hom(M, M'[n])$ form an effaceable $\delta$-functor, so they are universal and agree with $\Ext^n(M, M')$ for all $n \geq 0$. Indeed, for $n \leq 2$ the groups are effaceable since they agree with $\Ext$'s by the above, for $n > 2$ the groups are zero by \refsa{cohomdimAT}. 
\xpf

The functor $\eta_*: \DM(F) \r \DM(S)$ does not preserve Artin-Tate motives: 
$$\Hom_{\DM(S)}(\one, \eta_* \one(1)[1]) = \Hom_{\DM(F)}(\one, \one(1)[1]) = K^1(F)^{(1)}_\Q = F^\x \t \Q,$$ 
which is a countably infinite-dimensional $\Q$-vector space. However, the dimensions of all $\Hom$-groups in $T(S)$ are finite (\refle{morphisms}). This example is sharpened by the following proposition. It might be paraphrased by saying that the ``site'' of mixed Artin-Tate motives over $S$ has enough points. 

\satz \mylabel{satz_enoughpoints}
For any Artin-Tate motive $M$ over $S \subset \SpecOF$, the following are equivalent:
\begin{enumerate}[(i)]
\item \mylabel{item_M0} $M = 0$.
\item \mylabel{item_Meta} $M = \eta_* M_\eta$, where $M_\eta$ is some geometric motive over $F$.
\item \mylabel{item_Miupstar} $i_\pp^* M = 0$ for all closed points $\pp$ of $S$.
\item \mylabel{item_Miupshreak} $i_\pp^! M = 0$ for all closed points $\pp$ of $S$.
\end{enumerate}
\xsatz
\pf
The equivalence of \refit{Meta}, \refit{Miupstar}, and \refit{Miupshreak} is an easy consequence of Verdier duality on compact objects and the limiting localization triangle \refeq{localizationgeneric}, p.\ \pageref{eqn_localizationgeneric}. 
We now show \refit{Miupstar} $\Rightarrow$ \refit{M0}. Using localization, the claim for $M$ is implied by the one for $j^* M$ for any open immersion $j$. Therefore we may assume $M \in \tilde T(S)$. Using the $(-1)$-exactness of $i_\pp^*: \tilde T(S) \r T(\Fpp)$ we can even assume $M \in \tilde T^0(S)$. Given a short exact sequence in the abelian category $\tilde T^0(S)$
$$0 \r A \r M \r B \r 0$$
with $\eta_* \eta^* M = M$, it follows that $\eta_* \eta^* A = A$ and likewise for $B$. This is shown as follows: for all closed points $\pp \in S$, ${i_\pp}_* i_\pp^! M = 0$ implies $i_\pp^! B = i_\pp^! A[1]$, by the full faithfulness of ${i_\pp}_*$. The long exact $\pH^-$-sequence and the cohomological amplitude of $i_\pp^!$ (\refth{exactness}) shows $\pH^0 i_\pp^! B = \pH^1 i_\pp^! A$ and all other $\pH^* i_\pp^! B$, $\pH^* i_\pp^! A$ vanish. However, for any $B \in \tilde T^0(S)$, $i_\pp^! B$ is in cohomological degree $1$ (as opposed to the general range $[0, 1]$): this may be checked on generators of $\tilde T^0_a(S)$ for all $a$, where it follows directly from the definitions (see the proof of \refth{exactness}). Thus $\pH^0 i_\pp^! B = 0$, whence $i_\pp^! B = i_\pp^! A[1] = 0$ for all $\pp$. 

Thus the statement for $M$ is implied by the one for $A$ and $B$. By the characterization of $\tilde T^0(S)$ of \refsa{prep1}, we therefore only need to check the statement for generators of $\tilde T^0_{-2n}(S)$. 

We first do this in the case of Tate motives. Then $\tilde T^0_{-2n}(S)$ consists of direct sums of motives $G := \one(n)[1]$. In that case the claim is clear, since none of the (nonzero) generators $G$ satisfy $\eta_* \eta^* G = G$: we can twist it so that $n = 1$. Then $\H^0 (\eta_* \eta^* G)$ is infinite-dimensional, namely the group of units in some number field (tensored with $\Q$), but $\H^0(G)$ is the group of units in some ring of $S$-integers, which are of finite rank. 

In the case of Artin-Tate motives, the category $\tilde T^0_{-2n}(S)$ is generated by means of direct sums and summands by motives $G:= \phi_* \one(n)[1]$, $\phi: V \r S$ finite and \'etale. Actually, we may assume $\phi$ is Galois: by the same argument as in the proof of \refle{split}, after shrinking $S$ sufficiently, $\one_{V}$ is a direct summand of $\tilde \phi_* \one$ where $\tilde \phi: \tilde V \r V$ is the map corresponding to some normal closure of the function field extension $k(V) / k(S)$. Let $M$ be a summand of $G$ satisfying $\eta_* \eta^* M = M$. There is a map $f: S' \r S$ such that $f^* M$ is a Tate motive, \refle{split}. By base-change and the preceding step, we get $f^* M = 0$. The map $\End(M) \subset 
\End(G) \stackrel a \r \End(f^* G)$ factors over $\End (f^* M) = 0$, so we have to show that $a$ is injective. This is done with the same argument as in the proof of \refle{split}: we may shrink $S$ so that $f$ is \'etale. Since $\phi$ is Galois, we have
$$\End(G) = \Hom (\one, \phi^* \phi_* \one) = \Hom (\one, \one^{\oplus \deg \phi})$$
and
$$\End(f^* G) = \Hom (\one, \phi'^* \phi'_* \one) = \Hom (\one, \one^{\oplus \deg \phi'}),$$
where $\phi'$ is the pullback of $\phi$ along $f$. It is also Galois and $\deg \phi = \deg \phi'$. 
\xpf

\section{Weights} \mylabel{sect_weights}

This section develops a notion of weights on (mixed) Artin-Tate motives. We follow Bondarko and H\'ebert \cite{Bondarko:Weights, Hebert:Structures} for the definition of weights of Artin-Tate motives. That framework allows basic compatibility statements of weights, for example their behavior under functoriality. For mixed Artin-Tate motives, we show that this formalism can be used to produce a functorial and strict weight filtration. Again, this underlines the similarity of mixed motives with perverse sheaves. The latter type of results are strictly stronger than the ones obtained in \oc, and extend the ones of Levine and Wildeshaus concerning (Artin\nobreakdash-)Tate motives over fields \cite{Levine:TateMotives, Wildeshaus:ATM}. Briefly, we first relate the definition of weights to the one familiar from the theory of (perverse) $\ell$-adic sheaves (\refle{weightsBBD}). The technical key point is determining the subquotients of motives of the form $f_* \one[1]$, where $f$ is a finite flat map with regular domain (\refle{subquotients}). Via \refsa{subquotients}, this is the essential point in proving the strictness of the weight filtration in \refth{weights}.    

The following definition (in the more general situation of geometric motives) and \refle{weighttriangle} are due to Bondarko \cite[Th. 2.1.1, 2.2.1]{Bondarko:Weights} and H\'ebert \cite[3.2, 3.8]{Hebert:Structures}.

\defi \mylabel{defi_weights}
Let $\mathcal C$ be a triangulated category and $C$ a set of objects of $\mathcal C$. The category $\Extcl (C) \subset \mathcal C$ is defined to be the smallest full triangulated subcategory that contains $C$ and is stable under extensions, i.e., such that for any distinguished triangle
$$C \r X \r C'$$
with $C$ and $C'$ in $\Extcl (C)$, $X$ is also in $\Extcl (C)$.

Let $S$ be an (open or closed) subscheme of $\SpecOF$. Let $\Tw{0}(S)$ be the idempotent completion of the additive category (i.e., closed under direct summands and finite direct sums) generated by $f_* \one (a)[2a]$, where $f: S' \r S$ is a finite map such that $S'$ is regular (of dimension $\leq \dim S$) and $a \in \Z$ is arbitrary.

Put $\Tw{m} := \Tw{0}[m]$ and let $\Tw{\leq m}(S)$ be the idempotent completion of $\Extcl (\cup_{l \leq m} \Tw{l}(S))$ and define $\Tw{\geq m}(S)$ similarly.\footnote{Bondarko and H\'ebert use different notations. We follow H\'ebert here, thus our $\Tw{\leq m}$ would be $\Tw{\geq m}$ in Bondarko's notation.} We write $\Tw{>n}$ for $\Tw{\geq n+1}$ etc.  
\xdefi




\bem
For a map $f$ as in the definition $D(f_* f^* \one)=f_! f^! \one(1)[2] = f_* f^* \one(1)[2]$ by absolute purity and the regularity assumption on the domain of $f$. Thus $D(\Tw{\leq m}(S))=\Tw{\geq -m}(S)$. 
\xbem

\lemm \mylabel{lemm_weighttriangle}
\begin{enumerate}[(i)]
 \item \mylabel{item_weightvanishing} 
For any $M_{\leq m} \in \Tw{\leq m}(S)$, $M_{\geq m+1} \in \Tw{\geq m+1}(S)$,
$$\Hom(M_{\leq m}, M_{\geq m+1})=0.$$
\item
For any $M \in T(S)$ and any $m \in \Z$ there is a (non-unique) triangle
$$M_{\leq m} \r M \r M_{\geq m+1}$$
where $M_{\leq m} \in \Tw{\leq m}(S)$, $M_{\geq m+1} \in \Tw{\geq m+1}(S)$.
\item
For any map $f: S' \r S$ between regular schemes which is either a finite map, an open immersion or a closed immersion, the functors $f^*$ and $f_!$ preserve the $\Tw{\leq -}$-subcategories, and dually for $f_*$ and $f^!$.
\item \mylabel{item_weightlocalization}
Let $j$ and $i$ be an open and complementary closed embedding. Then $M$ is in $\Tw{\leq m}(S)$ iff $j^* M$ and $i^* M$ are in the corresponding $\Tw{\leq m}$-subcategories.
\end{enumerate} 
\xlemm

Recall that a mixed $\ell$-adic sheaf $\mathcal F$ on a curve $C / \Fq$ is said to be of weights $\leq m$ iff all pullbacks $i^* \mathcal F$ to all closed points of $C$ have this property \cite[5.1.9]{BBD}. We now relate this approach to weights to \refde{weights}. 

\lemm \mylabel{lemm_weightsBBD}
For any motive $M \in T(S)$, the following are equivalent:
\begin{enumerate}[(i)]
 \item \mylabel{item_eins} $M \in \Tw{\leq 0}(S)$ 
\item \mylabel{item_zwei}  $i^* M \in \Tw{\leq 0}(\Fpp)$ for any closed point $i: \SpecFpp \r S$.
\end{enumerate}
\xlemm

\pf
The implication \refit{eins} $\Rightarrow$ \refit{zwei} is trivial. Conversely, by \refle{weighttriangle}\refit{weightlocalization} and the localization triangle it is sufficient to show the implication \refit{zwei} $\Rightarrow$ \refit{eins} for $M \in \tilde T(S)$. There is a non-canonical isomorphism $i^*M \cong \oplus_n \pH^n (i^* M)[-n]$ by the semi-simplicity of Artin-Tate motives over finite fields. This in turn is (canonically) isomorphic to $\oplus_n i^* (\pH^n M) [-n]$ by the $(-1)$-exactness of $i^*$ (restricted to $\tilde T(S)$). Thus $i^* M \in \Tw{\leq 0}$ implies $i^* \pH^n M \in \Tw{\leq n}$ for all $n$. By definition, $\pH^n M \in \tilde T_{\wei{\leq n}}$ for all $n$ implies $M \in \tilde T_{\wei{\leq 0}}$. Thus we can replace $M$ by the $\pH^n M$ and show the statement for $M \in \tilde T^0(S)$ only. Using the same argument with respect to the auxiliary $t$-structure on $\tilde T^0(S)$ we reduce to $M \in \tilde T^0_{-2a}(S)$, $a \in \Z$. In this category, any extension splits. That is, its objects are direct summands of sums of motives $f_* \one(a)[1]$, where $f$ is finite etale. The implication \refit{zwei} $\Rightarrow$ \refit{eins} is obvious in that case.     
\xpf

The following lemma is the key stepstone in establishing the strictness of the weight filtration below. 

\lemm \mylabel{lemm_subquotients}
Let $f : \tilde S \r S$ and $g: V \r S$ be finite maps with regular domains $\tilde S$ and $V$ of dimension one and zero, respectively.  
Any subquotient in the abelian category $T^0(S)$ of $W := f_* \one[1]$ is a direct factor of $W$. In particular, it is also contained in $\Tw{\leq 1}$. A similar statement holds for $g_* \one$.
  
\xlemm 
\pf
Let $X \subset W$ be a subobject, $Y := W / X$. Let $j: S' \r S$ be an open immersion such that $j^* X$, $j^* W$, and $j^* Y \in \tilde T^0(S')$. Let $i : Z \r S$ be its complement. We have a commutative diagram in $T^0(S)$ with exact rows and columns:
$$\xymatrix
{
0 \ar[r] & i_* \pH^0 i^! X \ar[d] \ar[r] & i_* \pH^0 i^! W \ar[d] \\ 
0 \ar[r] & X \ar[r] \ar[d] & W \ar@/_/[l]_\tau \ar[r] \ar[d] & Y \ar[r] \ar[d] & 0 \\
0 \ar[r] & j_* j^* X \ar[d] \ar[r] & j_* j^* W \ar@/_/[l]_\sigma \ar[r] & j_* j^* Y \ar[r]  & 0 \\
& i_* \pH^1 i^! X. 
} 
$$
First of all, $i_* \pH^0 i^! W = 0$ by absolute purity, using that the domain of $f$ is regular. Thus $X \r j_* j^* X$ is a monomorphism. The curved arrows are splittings of the corresponding short exact sequences. Their existence is seen as follows: the third row exact sequence splits since $j^* X \in \tilde T^0_0(S')$ (which is a semi-simple category, since we use rational coefficients) and likewise for $j^* Y$, as follows from the $(-1)$-exactness of $p^*$ (restricted to $\tilde T^0(S')$), $(+1)$-exactness of $p^!$ for all closed points $p$ of $S'$ and purity. The map $W \r j_* j^* W \stackrel \sigma \r j_* j^* X \r i_* \pH^1 i^! X$ is zero: its image is a quotient of $W$ of the form $i_* N$ with $N \in T^0(Z)$. By the right-exactness of $i^*$, $N$ is a quotient of $\pH^0 i^* W = 0$. Thus $N = 0$. Hence $W \r j_* j^* W \stackrel \sigma \r j_* j^* X$ factors over some map $\tau: W \r X$. A short diagram chase shows that $\tau$ is a splitting of the second row exact sequence.   
%

The well-known second statement is easier. The details are omitted.
\xpf

The following corollary is a motivic analog of \cite[Th. 4.3.1 (i)]{BBD}.

\coro \mylabel{coro_ArtinianNoetherian} 
The category of mixed Artin-Tate motives is Artinian and Noetherian: given any $M \in T^0(S)$, and any sequence of subobjects in $T^0(S)$
$$0 = M_{-\infty} \subset \dots M_i \subset M_{i+1} \subset \dots M_\infty := M$$
there is an $n > 0$ such that $M_i = M_{i+1}$ for all $\infty > |i| > n$ and dually for quotients.
\xcoro
\pf
More generally, we claim that for any motive $N \in T(S)$, there is some number $l(N)$ such that the length of all subquotients of $\pH^n N$, for all $n \in \Z$, is bounded by $l(N)$. Given a triangle $N' \r N \r N''$ with $N'$ and $N''$ satisfying this claim, the claim also holds for $N$. This follows from the long exact cohomology sequence. Thus, it is sufficient to prove the claim for motives $f_* \one(a)[n]$, where $a, n \in \Z$ and $f$ is a finite map with regular domain (of dimension one or zero). For them, \refle{subquotients} and the compactness of $f_* \one$ settle the claim.  
\xpf

\satz \mylabel{satz_weightexact}
The truncation functors for the motivic $t$-structure are weight exact, that is to say
$$
\pH^s \Tw{\leq n}(S) = \Tw{\leq s+n}^0(S)
$$
and likewise for $\geq n$. In particular we have the following analog of \cite[5.1.8.]{BBD}: an Artin-Tate motive $M$ is of weights $\leq n$ iff all $\pH^s M$ are of weights $\leq s+n$. 
\xsatz
\pf
We can assume $s=0$. We clearly have $\Tw{\leq n}^0(S) = \pH^0 \Tw{\leq n}^0(S) \subseteq \pH^0 \Tw{\leq n}(S)$. Conversely, let $M \in \Tw{\leq n}(S)$. By \refle{weightsBBD} we have to show $i^* \pH^0 M$ is of weights $\leq n$ for all closed points $i$, if $i^* M$ has this property. The truncation functors $\tau_-$ of the motivic $t$-structure give a distinguished triangle 
$$i^* \tau_{\leq 0} M \r i^* M \r i^* \tau_{\geq 1} M \stackrel \delta \r i^* \tau_{\geq 0} M[1].$$
The third term is in cohomological degrees $0$ and $1$, the fourth one is in degrees $-2$ and $-1$. Hence the boundary map $\delta$ vanishes, by semi-simplicity of Artin-Tate motives over finite fields. Therefore the triangle splits (in a non-canonical way). The middle term $i^* M$ being of weights $\leq n$, the same follows for the summands $i^* \tau_{\leq 0} M$ and $i^* \tau_{\geq 1} M$. An induction shows that $i^* \pH^0 M$ is of weights $\leq n$.  

The concluding statement is a consequence of this and the truncation triangles for the motivic $t$-structure. 
\xpf

\satz \mylabel{satz_subquotients}
The category $\Tw{\leq n}^0(S)$ is stable under subquotients (that exist in the abelian category $T^0(S)$).
\xsatz
\pf
Recall that $\Tw{\leq n}^0(S) \stackrel{\text{\ref{satz_weightexact}}}= \pH^0 \Tw{\leq n}(S)$ is $\pH^0$ applied to the $\Ext$-closure (\refde{weights}) of the idempotent completion of 
$$\left \{ 
\begin{array}{c}
f_* \one(a)[2a+l]; \ \ l \leq n, \ a \in \Z, \\
f : S' \r S \text{ finite with regular domain}
\end{array}
\right \}.$$  
The subquotients of $\pH^0(f_* \one(a)[2a+l])$ are contained in $\Tw{\leq n}^0$ by \refle{subquotients}. It is thus sufficient to show the following statement: for any triangle
$$A \r X \r B$$ 
such that all subquotients of $\pH^0 A$ and $\pH^0 B$ are in $\Tw{\leq n}^0(S)$, all subquotients of $\pH^0 X$ are in $\Tw{\leq n}^0(S)$, too. Let $Y$ be a subobject of $\pH^0 X$. The triangle induces 
a short exact sequence
$$0 \r \coker (\pH^{-1} B \r \pH^0 A) \stackrel v \r \pH^0 X \stackrel w \r \ker (\pH^0 B \r \pH^{1} A) \r 0$$
which in turn induces
\eqn \mylabel{eqn_proof}
0 \r v^{-1} (Y) \r Y \r w (Y) \r 0.
\xeqn
The outer terms are subquotients of $\pH^0 A$ and $\pH^0 B$, respectively, hence they are in $\Tw{\leq n}^0(S)$. This category is stable under extensions in $T^0(S)$ by \refle{weightsBBD}. Therefore $Y \in \Tw{\leq n}^0(S)$. Quotients of $\pH^0 X$ are treated dually. 
%
\xpf

The following lemma is a partial converse to the general vanishings in weight structures (\refle{weighttriangle}\refit{weightvanishing}). It is used for the strictness of the weight filtration below.  

\lemm \mylabel{lemm_vanishconverse}
Let $n \in \Z$, and any $M \in \Tw{>0}^0(S)$, $M' \in \Tw{\leq 0}^0(S)$. Then
$$\Hom(M, M') = 0.$$
\xlemm
\pf
The argument in the above proof, cf.~\refeq{proof}, can be recycled to show the following: given any distinguished triangle $A \r B \r C$ such that for all subquotients $X$ of $\pH^0 A$ and of $\pH^0 C$ the group $\Hom(X, M')$ vanishes, the same vanishing holds for all subquotients $X$ of $\pH^0 B$, too. The dual statement of this, the description of $\Tw{-}^0(S) = \pH^0 \Tw{-}(S)$ in terms of extensions and the classification of subquotients in \refle{subquotients} show that it is suffices to see $\Hom(M, M')=0$, where 
$$M = \left \{ \begin{array}{c}
             f_* \one(a)[1] \ \  1-2a > 0 \\
\text{or} \\
    g_* \one(a), \ \ -2a > 0
\end{array} \right. 
\ \ \text {and} \ \ 
M' = \left \{ \begin{array}{c}
             f'_* \one(a')[1] \ \  1-2a' \leq 0 \\
\text{or} \\
    g'_* \one(a'), \ \ -2a' \leq 0
\end{array} \right.
$$ 
Here $f : V \r S$ and $f' : V' \r S$ are finite maps with regular domain of dimension $1$, $g$ and $g'$ are finite maps with $0$-dimensional image. Let $b = a' - a$. 
Most cases are an immediate consequence of absolute purity, except for the vanishing of $\Hom(f_* \one(a), f'_* \one(a'))$. It reduces to showing 
$$\Hom_{V'}(\tilde f_* \one, \one(b)) \left ( = \Hom_{W}(\one, \tilde f^! \one(b)) \right ) = 0,$$
where $\tilde f: W:= V \x_S V' \r V'$. If $W$ happens to be regular, this group identifies by absolute purity ($\tilde f^! \one = \tilde f^* \one$) with $K_{2b}(W)^{(b)}_\Q = 0$, since $b > 0$. In general there is the following argument, due to H\'ebert \cite[Theorem 3.1.]{Hebert:Structures} and Bondarko \cite[Lemma 1.1.4]{Bondarko:Weights}: let $n: W' \r W$ be the normalization map, $i: Z \subset W$ the exceptional ``divisor'', $Z'$ its preimage in $W'$, $z: Z' \r W$. The distinguished triangle $\one_W \r i_* \one_Z \oplus n_* \one_{W'} \r z_* \one_{Z'}$ induces an exact sequence
$$\dots \r \Hom(\tilde f_* i_* \one \oplus \tilde f_* n_* \one, \one(b)) \r \Hom(\tilde f_* \one, \one(b)) \r \Hom(\tilde f_* z_* \one[-1], \one(b)) \r \dots$$ 
The first half of the first term vanishes because of $b > 0$, the second one by the previous point. The last term vanishes for reasons of cohomological dimension.
\xpf

We can now construct the weight filtration. In a nutshell, the theorem says that weights for mixed Artin-Tate motives behave as they should, that is, as they do for mixed perverse $\ell$-adic sheaves \cite[5.3.5]{BBD} and mixed Hodge structures \cite[2.3.5]{Deligne:Hodge2}. The definition of the weight filtration $W_n M$ as the biggest subobject of $M$ of weight $\leq n$ is akin to a similar definition of Huber concerning the slice filtration of Hodge structures \cite[Lemma 2.1]{Huber:Slice}. 
It is worth noting that the classical proofs of the strictness of the weight filtration for perverse sheaves on a curve $C$ over $\Fq$ make use of the structural map $C \r \Spec \Fq$. In our situation, absolute purity (and the regularity of the base schemes we work over) take the r\^ole of this geometric piece of information. 

\theo \mylabel{theo_weights}
Let $M \in T^0(S)$ be a mixed Artin-Tate motive and $n \in \Z$. 
\begin{enumerate}[(i)]
\item \mylabel{item_a}
The set 
$$\mathcal W_n M := \{A \in \Tw{\leq n}^0, A \text{ is a subobject of }M \} / \text{isomorphism}$$
has a unique maximal element. Any choice of representatives of it is denoted $W_n (M)$. The assignment $M \mapsto W_n M$ is functorial (up to isomorphism, as $W_n M$ is only defined up to isomorphism).
\item \mylabel{item_d}
The quotient $W_{>n} M := M / W_n M$ is in $\Tw{>n}^0(S)$.  
\item \mylabel{item_c}
The \Def{weight filtration} $W_* -$ is strict: given any morphism $m: M \r M'$ in $T^0$, $\im W_n m = \im m \cap W_n M'$. Here $\im$ denotes the image of a map (in the abelian category $T^0(S)$).
\item \mylabel{item_e}
The assignment 
$$M \mapsto M := W_n M / W_{n-1} M$$
is an exact functor $\gr_n^W : T^0(S) \r \Tw{n}^0(S)$.
\end{enumerate}
\xtheo
\pf
\refit{a}: The existence of some maximal element in $\mathcal W_n M$ is assured by \refcor{ArtinianNoetherian}. 
Let $\iota_i : A_i \r M$, $i=1, 2$ be two maximal elements in $\mathcal W_n M$. Let $A := \im \left (\iota_1 \oplus \iota_2: A_1 \oplus A_2  \r M \right )$. This image is taken in the abelian category $T^0(S)$. As a quotient of $A_1 \oplus A_2$, $A$ is in $\Tw{\leq n}^0(S')$ by \refsa{subquotients}. Hence $A \in \mathcal W_n M$. By maximality of the $A_i$ we have $A_2 = A = A_1$.

Given a map $m: M_1 \r M_2$, fix representatives for $W_n M_i$, $i=1,2$.  
Again by \refsa{subquotients}, $W_n M_1 \subset M_1 \r M_2$ factors (uniquely, once representative are chosen) over $W_n M_2$. 
As the $W_n M_i$ are subobjects of $M_i$, the compatibility of $W_n$ with compositions is clear. 

\refit{d}: Let $A_{\wei{\leq n}} \r W_{>n} M \r A_{\wei{>n}}$ be any distinguished triangle where the outer terms are in $\Tw{\leq n}$ and $\Tw{>n}$, respectively (\refle{weighttriangle}). The induced long exact sequence 
$$\dots \stackrel{k}\r \pH^0 (A_{\wei{\leq n}}) \r W_{>n} M \r \pH^0 (A_{\wei{>n}}) \stackrel{l}\r \pH^1 (A_{\wei{\leq n}}) \r \dots$$
gives
$$0 \r \coker k \r W_{>n} M \r \ker l \r 0.$$
By \refsa{subquotients}, $V := \coker k \in \Tw{\leq n}$. Consider the pullback of the bottom row by $V$
$$\xymatrix{
0 \ar[r] &  
W_n M \ar[r] \ar@2{-}[d] &
V' \ar[r] \ar@{>->}[d]  &
V \ar[r] \ar@{>->}[d] &
0 \\
0 \ar[r] &
W_n M \ar[r] &
M \ar[r]&
W_{>n} M \ar[r] &
0. }
$$
As an extension of $V$ and $W_n M$, $V'$ is of weights $\leq n$, but also a subobject of $M$, so $V' = W_n M$. This shows $V = 0$, so that $W_{>n} M = \ker l \in \Tw{>n}$.

\refit{c}: 
We can assume $m$ is surjective and $M' = W_n M'$. We have to show $W_n M \r M'$ is surjective. 
Consider the commutative diagram with exact rows
$$\xymatrix{
0 \ar[r] &
W_n M \ar[r] \ar[d]^a &
M \ar[r]  \ar@{->>}[d]^m &
W_{>n} M \ar[r] \ar[d] &
0 \\
0 \ar[r] &
W_n M' \ar@2{-}[r] &
M' \ar[r] &
0.}
$$
The cokernel of $a$ is, as a quotient of $W_n M'$, in $\Tw{\leq n}^0(S)$. By \refit{d}, $W_{>n} M \in \Tw{>n}^0(S)$. \refle{vanishconverse} and the snake lemma imply $\coker a=0$. 

\refit{e}: By \refit{d} and the exact sequence 
\eqn \mylabel{eqn_weightexact}
0 \r W_{n-1} M \r W_n M  \r \gr_n^W M \r 0,
\xeqn
$\gr_n^W$ does map to $\Tw{n}^0(S)$. The exactness of $\gr_n^W$ is a well-known reformulation of the strictness of the weight filtration \cite[p.~8]{Deligne:Hodge2}.
\xpf

\bsp
Let $j: S' \r S$ be some open embedding with complement $i$. The sequence \refeq{weightexact} for $M := j_* j^* \one[1]$ and $n=2$ reads
$$0 \r \one[1] \r M \r i_* \one(-1) \r 0.$$ 
\xbsp

\bibliography{bib}  

\end{document}